\documentclass[11pt]{article}
\usepackage{amsmath,amssymb,latexsym,color,epsfig}
\usepackage[T1]{fontenc}
\begin{document}
\date{}
\title{The choosability version of Brooks' theorem --- a short proof}
\author{Michael Krivelevich\thanks{
School of Mathematical Sciences,
Tel Aviv University, Tel Aviv 6997801, Israel.
Email: {\tt krivelev@tauex.tau.ac.il}.  }
}
\vspace{-2cm}
\maketitle
\begin{abstract}
We present a short and self-contained proof of the choosability version of Brooks' theorem.
\end{abstract}
The following choosability version of Brooks' theorem is due to Vizing \cite{Viz76} and to Erd\H{o}s, Rubin and Taylor \cite{ERT79}.

\medskip

\noindent{\bf Theorem.}
Let $\Delta\ge 3$ be an integer, and let $G\ne K_{\Delta+1}$ be a connected graph of maximum degree at most $\Delta$. Then $G$ is $\Delta$-choosable.

\medskip
\noindent (The case $\Delta=2$ should -- and can easily -- be treated separately for this proof.)

\medskip

\noindent{\bf Proof.} The proof borrows its main idea from the nice argument of Zaj\k{a}c \cite{Zaj18} for the classical Brooks' theorem.
We proceed by induction on $n=|V(G)|$. The basic case $n\le \Delta$ is obvious --- given a list assignment $L$ for $V(G)$, one can just choose distinct colors for all vertices of $G$.

For the induction step, assume we are given a graph $G=(V,E)$ on $n$ vertices and a list assignment $L$ for the vertices of $G$ satisfying $|L(v)|=\Delta$ for all $v\in V$. We aim to find an $L$-coloring $f$ of $G$, which is a choice $f(v)\in L(v)$, $v\in V$, such that no edge of $G$ is monochromatic under $f$.

If $G$ contains a vertex $v$ with $d(v)<\Delta$, we can apply induction to color every connected component of $G-v$ from its lists in $L$. Let $f$ be the obtained coloring. It can be extended to $v$ by choosing $f(v)\in L(v)-\{f(u): (u,v)\in E\}$. We may thus assume $G$ is $\Delta$-regular.

Consider first the following special case: $G$ has a cycle $C$ on $k<n$ vertices and a vertex on the cycle having no neighbors outside $C$. Since $G$ is connected we can find two vertices $u,v$ adjacent along $C$ and such that $v$ has all its neighbors in $C$, and $u$ has some neighbor $w$ outside $C$. By the induction hypothesis the subgraph $G-V(C)$ is $\Delta$-choosable. Let $f$ be an $L$-coloring of $G-V(C)$. We now extend $f$ to $V(C)$. If $f(w)\not\in L(u)$, we choose $f(v)\in L(v)$ arbitrarily. If $f(w)\in L(u)\cap L(v)$, we set $f(v)=f(w)$. Finally, if $f(w)\in L(u)\setminus L(v)$, the lists $L(u)$ and $L(v)$ are different, and we can find $c\in L(v)\setminus L(u)$. We then set $f(v)=c$. In all three cases:
\begin{equation}\label{eq1}
|L(u)\cap \{f(v), f(w)\}|\le 1\,.
\end{equation}

Enumerate the vertices of $C$  when moving from $v$ to $u$ along $C$ as $V(C)=(v_1,\ldots,v_k)$ with $v_1=v$ and $v_k=u$.  Now we pick colors for $V(C)-\{v_1\}$ in the order $v_2,\ldots,v_k$. For $2\le i\le k-1$, vertex $v_i$ has at least one neighbor following it in this order, and thus an available color $f(v_i)$ can be found in its list $L(v_i)$. For $i=k$, at most $\Delta-1$ distinct colors from $L(v_k)$ have been used by $f$ on the neighbors of $v_k=u$ due to (\ref{eq1}), and we can choose $f(v_k)\in L(v_k)$ without creating a monochromatic edge.

Now we treat the general case. Since $G$ is connected and is not a clique, we can locate $v_1,v_2,v_3\in V$ such that $(v_1,v_2), (v_2,v_3)\in E$, but $(v_1,v_3)\not\in E$. Let $P=(v_1,v_2,v_3\ldots, v_{\ell})$ be a longest path in $G$ starting with $v_1,v_2,v_3$. All neighbors of $v_{\ell}$ reside on $P$. Let $v_i$ be the neighbor of $v_{\ell}$ farthest from it along $P$. The cycle $C=(v_i,v_{i+1},\ldots,v_{\ell})$ then contains all neighbors of $v_{\ell}$. If $C$ is not Hamiltonian then we are done by the special case considered above. We may thus assume $\ell=n$ and $i=1$. Let $v_j$ be a neighbor of $v_2$ different from $v_1,v_3$ (here we use the assumption $\Delta\ge 3$). Fix the following order  $\sigma$ on $V$: $\sigma=(v_1,v_3,\ldots,v_{j-1},v_n,v_{n-1},\ldots,v_j,v_2)$. As before, we can choose $f(v_1)\in L(v_1)$, $f(v_3)\in L(v_3)$ so that
\begin{equation}\label{eq2}
|L(v_2)\cap \{f(v_1), f(v_3)\}|\le 1\,.
\end{equation}
We now choose colors for the rest of $V$ in the order of $\sigma$. Every $v_i$ other than $v_2$ has at least one neighbor following it in $\sigma$, and thus we can set $f(v_i)\in L(v_i)$ without creating a monochromatic edge. Finally, when arriving to color $v_2$, we can allocate it a color from $L(v_2)$ distinct from the colors assigned to its neighbors, by (\ref{eq2}).\hfill$\Box$
\

\end{document}